\documentclass[psamsfonts]{amsart}
\usepackage{amssymb,eucal,latexsym,xy,graphics}
\xyoption{all}
\xyoption{ps}

\theoremstyle{plain}

\newtheorem{lemma}{Lemma}[section]
\newtheorem{theorem}[lemma]{Theorem}
\newtheorem{proposition}[lemma]{Proposition}
\newtheorem{corollary}[lemma]{Corollary}

\newtheorem*{stat}{\name}
\newcommand{\name}{testing}

\theoremstyle{definition}
\newtheorem{definition}[lemma]{Definition}
\newtheorem{example}[lemma]{Example}
\newtheorem{problem}{Problem}

\theoremstyle{remark}
\newtheorem{remark}[lemma]{Remark}

\newenvironment{all}[1]{\renewcommand{\name}{#1}\begin{stat}}
                         {\end{stat}}

\newcommand{\qedc}{{\qed}~{\rm Claim~{\theclaim}.}}

\numberwithin{equation}{section}

\newcommand{\Case}[1]{\medskip\noindent\textbf{\textit{Case~{#1}.}}}

\newcommand{\pup}[1]{\textup{(}{#1}\textup{)}}

\newcommand{\ba}{{\boldsymbol{a}}}
\newcommand{\bb}{{\boldsymbol{b}}}
\newcommand{\bc}{{\boldsymbol{c}}}
\newcommand{\bd}{{\boldsymbol{d}}}
\newcommand{\be}{{\boldsymbol{e}}}
\newcommand{\bu}{{\boldsymbol{u}}}
\newcommand{\bv}{{\boldsymbol{v}}}
\newcommand{\bw}{{\boldsymbol{w}}}
\newcommand{\xf}{{\boldsymbol{f}}}
\newcommand{\xg}{{\boldsymbol{g}}}
\newcommand{\xh}{{\boldsymbol{h}}}

\newcommand{\Bow}{\mathcal{D}_{\scriptscriptstyle{\bowtie}}}

\newcommand{\Pow}{\mathfrak{P}}

\newcommand{\jirr}{join-ir\-re\-duc\-i\-ble}

\newcommand{\into}{\hookrightarrow}
\newcommand{\onto}{\twoheadrightarrow}

\newcommand{\set}[1]{\{#1\}}
\newcommand{\setm}[2]{\set{#1\mid#2}}
\newcommand{\seq}[1]{\langle{#1}\rangle}
\newcommand{\seqm}[2]{\seq{#1\mid#2}}
\newcommand{\famm}[2]{(#1\mid#2)}

\newcommand{\op}{\mathrm{op}}
\newcommand{\ol}[1]{\overline{#1}}
\newcommand{\oll}[1]{\,\overline{\!#1}}

\newcommand{\zero}{\mathbf{0}}
\newcommand{\one}{\mathbf{1}}
\newcommand{\two}{\mathbf{2}}
\newcommand{\go}{\omega}
\newcommand{\gos}{\go\setminus\set{0}}

\DeclareMathOperator{\Id}{Id}

\DeclareMathOperator{\J}{J}

\newcommand{\B}{\mathbf{B}}
\newcommand{\cA}{\mathcal{A}}

\newcommand{\cC}{\mathcal{C}}
\newcommand{\Dac}{\mathcal{D}_{\mathrm{ac}}}
\newcommand{\cF}{\mathcal{F}}
\newcommand{\cI}{\mathcal{I}}
\newcommand{\cV}{\mathcal{V}}
\newcommand{\cD}{\mathcal{D}}
\newcommand{\cS}{\mathcal{S}}
\newcommand{\cH}{\mathcal{H}}

\newcommand{\ZZ}{\mathbb{Z}}

\DeclareMathOperator{\Subc}{Sub_c}
\DeclareMathOperator{\Sub}{Sub}
\DeclareMathOperator{\NSubc}{NSub_c}
\DeclareMathOperator{\NSub}{NSub}
\DeclareMathOperator{\Conc}{Con_c}
\DeclareMathOperator{\Con}{Con}
\DeclareMathOperator{\Idc}{Id_c}
\newcommand{\Idl}{\operatorname{Id}^{\ell}}
\newcommand{\Idlc}{\operatorname{Id}^{\ell}_{\mathrm{c}}}

\newcommand{\jz}{$\langle\vee,0\rangle$}
\newcommand{\jzu}{$\langle\vee,0,1\rangle$}
\newcommand{\jzs}{\jz-semi\-lat\-tice}
\newcommand{\jzus}{\jzu-semi\-lat\-tice}
\newcommand{\jzh}{\jz-ho\-mo\-mor\-phism}

\newcommand{\jze}{\jz-em\-bed\-ding}
\newcommand{\jzue}{\jzu-em\-bed\-ding}

\newcommand{\jh}{join-ho\-mo\-mor\-phism}

\newcommand{\fin}[1]{[#1]^{<\go}}

\DeclareMathOperator{\Eq}{Eq}

\newcommand{\wurpe}{\mathrm{WURP}^=}

\begin{document}

\title[Congruence-permutable algebras]{Distributive congruence lattices of
congruence-permutable algebras}
 \author[P. R\r{u}\v{z}i\v{c}ka]{Pavel R\r{u}\v{z}i\v{c}ka}
 \author[J.~T\r{u}ma]{Ji\v r\'\i~T\r{u}ma}
 \address{Charles University in Prague\\
          Faculty of Mathematics and Physics\\
          Department of Algebra\\
          Sokolovsk\'a 83\\
          Charles University\\
          186 00 Praha 8\\
          Czech Republic}
 \email[P. R\r{u}\v{z}i\v{c}ka]{ruzicka@karlin.mff.cuni.cz}
 \urladdr[P. R\r{u}\v{z}i\v{c}ka]{http://www.karlin.mff.cuni.cz/\~{}ruzicka/}
 \email[J. T\r{u}ma]{tuma@karlin.mff.cuni.cz}
 \urladdr[J. T\r{u}ma]{http://www.karlin.mff.cuni.cz/\~{}tuma/}

\author[F.~Wehrung]{Friedrich~Wehrung}
\address{LMNO, CNRS UMR 6139\\
         D\'epartement de Math\'ematiques, BP 5186\\
         Universit\'e de Caen, Campus 2\\
         14032 Caen cedex\\
         France}
 \email{wehrung@math.unicaen.fr}
 \urladdr{http://www.math.unicaen.fr/\~{}wehrung}

 \date{\today}
\thanks{The first two authors were partially supported by the
institutional grant MSM 0021620839 and INTAS project 03-51-4110. The first
author was partially supported by grant GAUK 284/2003/B MAT/MFF and
post-doctoral grant GA\v CR 201/03/P140. The second author was partially
supported by grants GAUK 448/2004/B-MAT and GA\v CR 201/03/0937}

\begin{abstract}
We prove that every distributive algebraic lattice with at
most~$\aleph_1$ compact elements is isomorphic to the normal subgroup
lattice of some group and to the submodule lattice of some right module.
The $\aleph_1$ bound is optimal, as we find a distributive algebraic
lattice $D$ with~$\aleph_2$ compact elements that is not isomorphic to the
congruence lattice of any algebra with almost permutable congruences
(hence neither of any group nor of any module), thus solving negatively a
problem of E.\,T. Schmidt from~1969. Furthermore, $D$ may be taken as the
congruence lattice of the free bounded lattice on~$\aleph_2$ generators
in any non-distributive lattice variety.

Some of our results are obtained \emph{via} a functorial approach of the
semilattice-valued `distances' used by B. J\'onsson in his proof of
Whitman's embedding Theorem. In particular, the semilattice of compact
elements of $D$ is not the range of any distance satisfying the
V-condition of type~$3/2$. On the other hand, every distributive \jzs\ is
the range of a distance satisfying the V-condition of type~2.
This can be done \emph{via} a functorial construction.
\end{abstract}

\maketitle

\section*{Introduction}\label{S:Intro}

Representing algebraic lattices as congruence lattices of
algebras often gives rise to very hard open problems. The most well-known of those problems, the \emph{Congruence Lattice Problem}, usually abbreviated CLP, asks whether every distributive algebraic lattice is isomorphic to the congruence
lattice of some lattice, see the survey paper \cite{CLPSurv}. This problem has been solved recently by the third author in~\cite{CLP}.
For algebraic lattices that are not necessarily distributive, there are
several deep results, one of the most remarkable, due to W.\,A.
Lampe~\cite{Lamp82}, stating that \emph{every algebraic lattice with
compact unit is isomorphic to the congruence lattice of some groupoid}.
This result is further extended to join-complete, unit-preserving,
compactness preserving maps between two algebraic lattices \cite{Lamp05}.

Although some of our methods are formally related to Lampe's, for example
the proof of Theorem~\ref{T:ReprType2} \emph{via}
Proposition~\ref{P:ReprType2}, we shall be concerned only about
\emph{distributive} algebraic lattices. This topic contains some not so
well-known but also unsolved problems, as, for example, whether every
distributive algebraic lattice is isomorphic to the congruence lattice of
an algebra in some \emph{congruence-distributive variety}.

If one drops congruence-distributivity, then one would expect the
problems to become easier. Consider, for example, the two following
problems:

\begin{all}{CGP}
Is every distributive algebraic lattice isomorphic to the normal subgroup
lattice of some group?
\end{all}

\begin{all}{CMP}
Is every distributive algebraic lattice isomorphic to the submodule
lattice of some module?
\end{all}

The problem CGP was originally posed for \emph{finite} distributive
(semi)lattices by E.\,T. Schmidt as \cite[Problem~5]{SCSy}. A positive
solution was provided by H.\,L. Silcock, who proved in particular that
every finite distributive lattice~$D$ is isomorphic to the normal
subgroup lattice of some
\emph{finite} group~$G$ (see~\cite{Silc77}). P.\,P. P\'alfy proved later
that~$G$ may be taken
\emph{finite solvable} (see~\cite{Palf86}). However, the general question
seemed open until now. Similarly, the statement of CMP has been
communicated to the authors by Jan Trlifaj, and nothing seemed to be
known about the general case.

A common feature of the varieties of all groups and of all modules over a
given ring is that they are \emph{congruence-permutable}, for example,
any two congruences of a group are permutable. Thus both CGP and CMP
are, in some sense, particular instances of the following question:

\begin{all}{CPP (see \cite[Problem~3]{SCSy})}
Is every distributive algebraic lattice isomorphic to the congruence
lattice of some algebra with permuting congruences?
\end{all}

Although the exact formulation of \cite[Problem~3]{SCSy} asked whether
every \emph{Arguesian} algebraic lattice is isomorphic to the congruence
lattice of an algebra with permutable congruences, it was mentioned there
that even the
\emph{distributive} case was open. Meanwhile, the Arguesian case was
solved negatively by M.\,D. Haiman \cite{Haim87,Haim91}, however, the
distributive case remained open.

Recall that an algebra $A$ has \emph{almost permutable
congruences} (see \cite{TuWe}), if $\ba\vee\bb=\ba\bb\cup\bb\ba$, for all congruences $\ba$, $\bb\in\Con A$ (where the notation $\ba\bb$ stands for the usual composition of relations).
The three-element chain is an easy example of a lattice with almost permutable congruences but not with permutable congruences. On the other hand, it is not difficult to verify that \emph{every almost congruence-permutable variety of algebras is congruence-permutable}. The last two authors of the present
paper obtained in \cite{TuWe} negative congruence representation results of distributive semilattices by \emph{lattices} with almost permutable congruences, but nothing was
said there about arbitrary \emph{algebras} with permutable congruences.
Furthermore, our attempts based on the ``uniform refinement properties''
introduced in that paper failed, as these properties turned out to be
quite lattice-specific.

In the present paper, we introduce a general framework that makes it
possible to extend the methods of \cite{TuWe} to arbitrary algebras, and
thus solving CPP---and, in fact, its generalization to algebras
with almost permutable congruences---negatively. Hence, both CGP and CMP
also have negative solutions. In fact, the negative solution obtained in
CGP for \emph{groups} extends to \emph{loops}, as the variety of all
loops is also congruence-permutable. Another byproduct is that we also
get a negative solution for the corresponding problem for
\emph{lattice-ordered groups}, see also Problem~\ref{Pb:lGroups}.

Our counterexample is the same
as in~\cite{PTW} and in~\cite{TuWe}, namely the congruence lattice of a
free lattice with at least $\aleph_2$ generators in any non-distributive
variety of lattices. We also show that the size $\aleph_2$ is optimal, by
showing that \emph{every distributive algebraic lattice with at
most~$\aleph_1$ compact elements is isomorphic to the submodule lattice
of some module, and also to the normal subgroup lattice of some locally
finite group}, see Theorems~\ref{T:SubMAl1} and~\ref{T:Al1GrpRepr}. We
also prove that
\emph{every distributive algebraic lattice with at most countably many
compact elements is isomorphic to the $\ell$-ideal lattice of some
lattice-ordered group}, see Theorem~\ref{T:Reprlgroups}.

In order to reach our negative results, the main ideas are the following.
\begin{itemize}
\item[(1)] Forget about the algebraic structure, just keep the partition
lattice representation.

\item[(2)] State a weaker ``uniform refinement property'' that settles
the negative result.
\end{itemize}

For Point~(1), we are looking for a very special sort of lattice
homomorphism of a given lattice into some partition lattice, namely, the
sort that is induced, as in Proposition~\ref{P:Dist2Part}, by a
\emph{semilattice-valued distance}, see Definition~\ref{D:Vdist}. For a
\jzs\ $S$ and a set $X$, an~$S$-valued distance on $X$ is a map
$\delta\colon X\times X\to S$ satisfying the three usual statements
characterizing distances (see Definition~\ref{D:Vdist}). Every such
$\delta$ induces a map $\varphi$ from $S$ to the partition lattice of $X$
(see Proposition~\ref{P:Dist2Part}), and if $\delta$ satisfies the
so-called \emph{V-condition}, then $\varphi$ is a \jh. Furthermore, the
V-condition of type~$n$ says that the equivalences in the range of~$\varphi$ are pairwise $(n+1)$-permutable.
Those ``distances'' have been introduced by B. J\'onsson for providing a
simple proof of Whitman's Theorem that every lattice can be embedded into
some partition lattice, see \cite{Jons53} or Theorems~IV.4.4 and~IV.4.8
in~\cite{GLT2}.

While it is difficult to find a suitable notion of morphism
between partition lattices, it is easy to do such a thing with our
distances, see Definition~\ref{D:Vdist}. This makes it possible to define
what it means for a commutative diagram of \jzs s to have a lifting,
modulo the forgetful functor, by distances. In particular, we prove, in
Theorem~\ref{T:NoLifttyp1}, that the cube $\Dac$ considered in
\cite[Section~7]{TuWe} does not have a lifting by any diagram of
V-distances ``of type~$3/2$'', that is, the equivalences in the ranges of
the corresponding partition lattice representations cannot all be almost
permutable. This result had been obtained only for lattices
in~\cite{TuWe}.

The original proof of Theorem~\ref{T:NoLifttyp1} was our main inspiration
for getting a weaker ``uniform refinement property'', that we denote here
by $\wurpe$ (see Definition~\ref{D:VWURPe}). First, we prove that
if $\delta\colon X\times X\to S$ is an~$S$-valued V-distance of type~$3/2$
with range generating $S$, then~$S$ satisfies $\wurpe$ (see
Theorem~\ref{T:type152WURP}). Next, we prove that for any free
lattice~$F$ with at least~$\aleph_2$ generators in any non-distributive
variety of lattices, the compact congruence semilattice $\Conc F$ does not
satisfy $\wurpe$ (see Corollary~\ref{C:FreeNon2}). Therefore, $\Con F$ is
not isomorphic to $\Con A$, for any algebra~$A$ with almost permutable
congruences (see Corollary~\ref{C:FreeNon}).

On the positive side, we explain why all previous attempts at finding
similar negative results for representations of type~$2$ (and above)
failed. We prove, in particular, that \emph{for every distributive \jzs\
$S$, there exists a surjective V-distance
$\delta_S\colon X_S\times X_S\onto S$ of type~$2$, which, moreover,
depends functorially on~$S$} (see Theorem~\ref{T:ReprType2}). In
particular, the diagram $\Bow$ considered in \cite{Bowtie}, which is not
liftable, with respect to the congruence lattice functor, in any variety
whose congruence lattices satisfy a nontrivial identity, is nevertheless
liftable by V-distances of type~$2$.

\section*{Basic concepts}\label{S:Basic}

For elements $x$ and $y$ in an algebra $A$, we denote by $\Theta_A(x,y)$,
or $\Theta(x,y)$ if $A$ is understood, the least congruence of $A$
that identifies $x$ and $y$. Furthermore, in case~$A$ is a lattice, we put
$\Theta^+_A(x,y)=\Theta_A(x\wedge y,x)$.
We denote by $\Con A$ (resp., $\Conc A$) the lattice (resp., semilattice)
of all compact (i.e., finitely generated) congruences of $A$.

For join-semilattices $S$ and $T$, a join-homomorphism
$\mu\colon S\to T$ is \emph{weakly distributive} (see \cite{UnifRef}),
if for every $c\in S$ and $a$, $b\in T$, if $\mu(c)\leq a\vee b$, then
there are $x$, $y\in S$ such that $c\leq x\vee y$, $\mu(x)\leq a$, and
$\mu(y)\leq b$.

A \emph{diagram} in a category $\cC$ is a functor
$\mathbf{D}\colon\cI\to\cC$, for some category $\cI$. For a functor
$\mathbf{F}\colon\cA\to\cC$, a \emph{lifting} of $\mathbf{D}$ with
respect to $\mathbf{F}$ is a functor $\Phi\colon\cI\to\cA$ such that the
composition $\mathbf{F}\circ\Phi$ is naturally equivalent to $\mathbf{D}$.

For a set $X$ and a natural number $n$, we denote by $[X]^n$ the set
of all $n$-elements subsets of $X$, and we put
$[X]^{<\go}=\bigcup\famm{[X]^n}{n<\go}$. The following statement of
infinite combinatorics can be found in C. Kuratowski~\cite{Kura51}.

\begin{all}{The Kuratowski Free Set Theorem}
Let $n$ be a positive integer and let $X$ be a set. Then
$|X|\geq\aleph_n$ if{f} for every map $\Phi\colon[X]^n\to\fin{X}$,
there exists $U\in[X]^{n+1}$ such that 
$u\notin\Phi(U\setminus\set{u})$, for any
$u\in U$.
\end{all}

As in \cite{PTW,UnifRef}, only the case $n=2$ will be used.

We identify every natural number $n$ with the set
$\set{0,1,\dots,n-1}$, and we denote by~$\go$ the set of all
natural numbers.

\section{V-distances of type $n$}\label{S:Vdist}

\begin{definition}\label{D:Vdist}
Let $S$ be a \jzs\ and let $X$ be a set. A map
$\delta\colon X\times\nobreak X\to\nobreak S$ is an \emph{$S$-valued
distance} on $X$, if the following statements hold:
 \begin{enumerate}
 \item $\delta(x,x)=0$, for all $x\in X$.

 \item $\delta(x,y)=\delta(y,x)$, for all $x$, $y\in X$.

 \item $\delta(x,z)\leq\delta(x,y)\vee\delta(y,z)$, for all $x$, $y$,
 $z\in X$.
 \end{enumerate}
The \emph{kernel} of $\delta$ is defined as
$\setm{\seq{x,y}\in X\times X}{\delta(x,y)=0}$. The \emph{V-condition}
on $\delta$ is the following condition:
 \begin{quote}
For all $x$, $y\in X$ and all $\ba$, $\bb\in S$ such that
$\delta(x,y)\leq\ba\vee\bb$, there are $n\in\gos$ and $z_0=x$, $z_1$,
\dots, $z_{n+1}=y$ such that for all $i\leq n$,
$\delta(z_i,z_{i+1})\leq\ba$ in case $i$ is even, while
$\delta(z_i,z_{i+1})\leq\bb$ in case $i$ is odd.
 \end{quote}
In case $n$ is the same for all $x$, $y$, $\ba$, $\bb$, we say that the
distance $\delta$ satisfies the \emph{V-condition of type~$n$}, or is a
\emph{V-distance of type~$n$}.

We say that $\delta$ satisfies the \emph{V-condition of
type~$3/2$}, or is a \emph{V-distance of type~$3/2$}, if for all $x$,
$y\in X$ and all $\ba$, $\bb\in S$ such that
$\delta(x,y)\leq\ba\vee\bb$, there exists $z\in X$ such that either
($\delta(x,z)\leq\ba$ and $\delta(z,y)\leq\bb$) or ($\delta(x,z)\leq\bb$
and $\delta(z,y)\leq\ba$).

We say that a \emph{morphism} from $\lambda\colon X\times X\to A$ to
$\mu\colon Y\times Y\to B$ is a pair $\seq{f,\xf}$, where
$\xf\colon A\to B$ is a \jzh\ and $f\colon X\to Y$ is a map such that
$\xf(\lambda(x,y))=\mu(f(x),f(y))$, for all $x$, $y\in X$. The
\emph{forgetful functor} sends $\lambda\colon X\times X\to A$ to $A$ and
$\seq{f,\xf}$ to $\xf$.
\end{definition}

Denote by $\Eq X$ the lattice of all equivalence relations on a set $X$.
For a positive integer $n$, we say as usual that $\alpha$, $\beta\in\Eq X$
are \emph{$(n+1)$-permutable}, if
$\gamma_0\gamma_1\cdots\gamma_n=\gamma_1\gamma_2\cdots\gamma_{n+1}$,
where $\gamma_k$ is defined as~$\alpha$ if~$k$ is even and as~$\beta$
if~$k$ is odd, for every natural number~$k$. In particular, $2$-permutable is the same as permutable. With every distance is
associated a homomorphism to some $\Eq X$, as follows.

\begin{proposition}\label{P:Dist2Part}
Let $S$ be a \jzs\ and let $\delta\colon X\times X\to S$ be an~$S$-valued
distance. Then one can define a map $\varphi\colon S\to\Eq X$ by the rule
 \[
 \varphi(\ba)=\setm{\seq{x,y}\in X\times X}{\delta(x,y)\leq\ba},
 \text{ for all }\ba\in S.
 \]
Furthermore,
\begin{enumerate}
\item The map $\varphi$ preserves all existing meets.

\item If $\delta$ satisfies the V-condition, then $\varphi$ is a \jh.

\item If the range of $\delta$ join-generates $S$, then $\varphi$ is an
order-embedding.

\item If the distance $\delta$ satisfies the V-condition of type~$n$, then all equivalences in the range of $\varphi$ are pairwise $(n+1)$-permutable.
\end{enumerate}
\end{proposition}

Any algebra gives rise to a natural distance, namely the map
$\seq{x,y}\mapsto\Theta(x,y)$ giving the principal congruences.

\begin{proposition}\label{P:PermCong}
Let $n$ be a positive integer and let $A$ be an algebra with
$(n+1)$-permutable congruences. Then the semilattice
$\Conc A$ of compact congruences of $A$ is join-generated by the range
of a V-distance of type~$n$.
\end{proposition}

\begin{proof}
Let $\delta\colon A\times A\to\Conc A$ be defined by
$\delta(x,y)=\Theta_A(x,y)$, the principal congruence generated by
$\seq{x,y}$, for all $x$, $y\in A$. The assumption that $A$ has
$(n+1)$-permutable congruences means exactly that $\delta$ is a V-distance of
type~$n$.
\end{proof}

Of course, $A$ has almost permutable congruences if and only if the
canonical distance $\Theta_A\colon A\times A\to\Conc A$ satisfies the
V-condition of type~$3/2$.

We shall focus attention on three often encountered varieties all members
of which have permutable (i.e., $2$-permutable) congruences:
\begin{itemize}
\item[---] The variety of all right modules over a given ring $R$. The
congruence lattice of a right module $M$ is canonically isomorphic to the
submodule lattice $\Sub M$ of $M$. We shall denote by $\Subc M$ the \jzs\
of all finitely generated submodules of $M$.

\item[---] The variety of all groups. The
congruence lattice of a group $G$ is canonically isomorphic to the
normal subgroup lattice $\NSub G$ of $G$. We shall denote by $\NSubc G$
the \jzs\ of all finitely generated normal subgroups of $G$.

\item[---] The variety of all $\ell$-groups (i.e., lattice-ordered
groups), see \cite{AnFe}. The congruence lattice of an $\ell$-group $G$ is
canonically isomorphic to the lattice $\Idl G$ of all convex
normal subgroups, or \emph{$\ell$-ideals}, of $G$. We shall denote by
$\Idlc G$ the \jzs\ of all finitely generated $\ell$-ideals of $G$.
\end{itemize}

Hence we obtain immediately the following result.

\begin{corollary}\label{C:GrpMod}\hfill
\begin{enumerate}
\item Let $M$ be a right module over any ring $R$. Then $\Subc M$ is
join-generated by the range of a V-distance of type~$1$ on $M$.

\item Let $G$ be a group. Then $\NSubc G$ is join-generated by
the range of a V-distance of type~$1$ on $G$.

\item Let $G$ be an $\ell$-group. Then $\Idlc G$ is join-generated by
the range of a V-distance of type~$1$ on $G$.
\end{enumerate}
\end{corollary}

The V-distances corresponding to (i), (ii), and (iii) above are,
respectively, given by $\delta(x,y)=(x-y)R$, $\delta(x,y)=[xy^{-1}]$
(the normal subgroup of $G$ generated by $xy^{-1}$), and
$\delta(x,y)=G(xy^{-1})$ (the $\ell$-ideal of $G$ generated by $xy^{-1}$).

The assignments $M\mapsto\Subc M$, $G\mapsto\NSubc G$, and
$G\mapsto\Idlc G$ can be canonically extended to direct limits preserving
functors to the category of all \jzs s with \jzh s.

\section{An even weaker uniform refinement property}\label{S:WURPe}

The following infinitary axiom $\wurpe$ is a weakening of
all the various ``uniform refinement properties'' considered in
\cite{PTW,TuWe,UnifRef}. Furthermore, the proof that follows, aimed at
obtaining Theorem~\ref{T:NonUnif}, is very similar to the proofs of
\cite[Theorem~3.3]{PTW} and \cite[Theorem~2.1]{TuWe}.

\begin{definition}\label{D:VWURPe}
Let $\be$ be an element in a \jzs\ $S$. We say that~$S$ satisfies
$\wurpe(\be)$, if there exists a positive integer $m$ such that
for all families $\seqm{\ba_i}{i\in I}$ and
$\seqm{\bb_i}{i\in I}$ of elements of~$S$ such that
$\be\leq\ba_i\vee\bb_i$ for all $i\in I$, there are a
$m$-sequence $\seqm{I_u}{u<m}$ of subsets of $I$ such that
$\bigcup\famm{I_u}{u<m}=I$ and a family
$\seqm{\bc_{i,j}}{\seq{i,j}\in I\times I}$ of elements of $S$ such that
the following statements hold:
 \begin{enumerate}
 \item $\bc_{i,j}\leq\ba_i\vee\ba_j$ and $\bc_{i,j}\leq\bb_i\vee\bb_j$,
 for all $u<m$ and all $i$, $j\in I_u$.

 \item $\be\leq\ba_j\vee\bb_i\vee\bc_{i,j}$,
 for all $u<m$ and all $i$, $j\in I_u$.

 \item $\bc_{i,k}\leq\bc_{i,j}\vee\bc_{j,k}$,  for all $i$, $j$,
 $k\in I$.
 \end{enumerate}
Say that $S$ satisfies $\wurpe$, if $S$ satisfies $\wurpe(\be)$ for
all $\be\in S$.
\end{definition}

The following easy lemma is instrumental in the proof of
Corollary~\ref{C:FreeNon}.

\begin{lemma}\label{L:WDUnif}
Let $S$ and $T$ be \jzs s, let $\mu\colon S\to T$ be a weakly
distributive \jzh, and let $\be\in S$. If $S$ satisfies
$\wurpe(\be)$, then $T$ satisfies $\wurpe(\mu(\be))$.
\end{lemma}

\begin{theorem}\label{T:type152WURP}
Let $S$ be a \jzs\ and let $\delta\colon X\times X\to S$ be
a V-distance of type~$3/2$ with range join-generating $S$. Then $S$
satisfies $\wurpe$.
\end{theorem}

\begin{proof}
Let $\be\in S$.
As $S$ is join-generated by the range of $\delta$, there are a positive
integer $n$ and elements $x_\ell$, $y_\ell\in X$, for $\ell<n$, such
that $\be=\bigvee\famm{\delta(x_\ell,y_\ell)}{\ell<n}$. For all $i\in I$
and all $\ell<n$, from $\delta(x_\ell,y_\ell)\leq\ba_i\vee\bb_i$ and the
assumption on $\delta$ it follows that there exists $z_{i,\ell}\in X$
such that
 \begin{equation}\label{Eq:APCdel}
 \begin{aligned}
 \text{either }\delta(x_\ell,z_{i,\ell})\leq\ba_i&\text{ and }
 \delta(z_{i,\ell},y_\ell)\leq\bb_i\\
 \text{or }\delta(x_\ell,z_{i,\ell})\leq\bb_i&\text{ and }
 \delta(z_{i,\ell},y_\ell)\leq\ba_i.
 \end{aligned}
 \end{equation}
For all $i\in I$ and all $\ell<n$, denote by $P(i,\ell)$ and
$Q(i,\ell)$ the following statements:
 \[
 \begin{aligned}
 P(i,\ell):\quad&\delta(x_\ell,z_{i,\ell})\leq\ba_i\text{ and }
 \delta(z_{i,\ell},y_\ell)\leq\bb_i;\\
 Q(i,\ell):\quad&\delta(x_\ell,z_{i,\ell})\leq\bb_i\text{ and }
 \delta(z_{i,\ell},y_\ell)\leq\ba_i.
 \end{aligned}
 \]
We shall prove that $m=2^n$ is a suitable choice for witnessing
$\wurpe(\be)$. So let~$U$ denote the powerset of $n$, and put
 \[
 I_u=\setm{i\in I}{(\forall\ell\in u)P(i,\ell)\text{ and }
 (\forall\ell\in n\setminus u)Q(i,\ell)},\text{ for all }u\in U.
 \]
We claim that $I=\bigcup\famm{I_u}{u\in U}$. Indeed, let $i\in I$,
and put $u=\setm{\ell<n}{P(i,\ell)}$. It follows from \eqref{Eq:APCdel}
that $Q(i,\ell)$ holds for all $\ell\in n\setminus u$, whence $i\in I_u$.
Now we put
 \[
 \bc_{i,j}=\bigvee\famm{\delta(z_{i,\ell},z_{j,\ell})}{\ell<n},
 \text{ for all }i,\,j\in I,
 \]
and we prove that the family $\seqm{\bc_{i,j}}{\seq{i,j}\in I\times I}$
satisfies the required conditions, with respect to the family
$\seqm{I_u}{u\in U}$ of $2^n$ subsets of $I$.
So, let $i$, $j$, $k\in I$. The inequality
$\bc_{i,k}\leq\bc_{i,j}\vee\bc_{j,k}$ holds trivially.

Now suppose that $i$, $j\in I_u$, for some $u\in U$.

Let $\ell<n$. If $\ell\in u$, then
 \begin{align*}
 \delta(z_{i,\ell},z_{j,\ell})&\leq
 \delta(z_{i,\ell},x_\ell)\vee\delta(x_\ell,z_{j,\ell})\leq
 \ba_i\vee\ba_j,\\
 \delta(x_\ell,y_\ell)&\leq
 \delta(x_\ell,z_{j,\ell})\vee\delta(z_{j,\ell},z_{i,\ell})\vee
 \delta(z_{i,\ell},y_\ell)\leq\ba_j\vee\bc_{i,j}\vee\bb_i,
 \end{align*}
while if $\ell\in n\setminus u$, 
 \begin{align*}
 \delta(z_{i,\ell},z_{j,\ell})&\leq
 \delta(z_{i,\ell},y_\ell)\vee\delta(y_\ell,z_{j,\ell})
 \leq\ba_i\vee\ba_j,\\
 \delta(x_\ell,y_\ell)&\leq
 \delta(x_\ell,z_{i,\ell})\vee\delta(z_{i,\ell},z_{j,\ell})\vee
 \delta(z_{j,\ell},y_\ell)\leq\bb_i\vee\bc_{i,j}\vee\ba_j.
 \end{align*}
whence both inequalities $\delta(z_{i,\ell},z_{j,\ell})\leq\ba_i\vee\ba_j$
and  $\delta(x_\ell,y_\ell)\leq\ba_j\vee\bb_i\vee\bc_{i,j}$ hold in any
case. Hence $\bc_{i,j}\leq\ba_i\vee\ba_j$ and
$\be\leq\ba_j\vee\bb_i\vee\bc_{i,j}$. Exchanging~$x$ and~$y$ in the
argument leading to the first inequality also yields that
$\bc_{i,j}\leq\bb_i\vee\bb_j$.
\end{proof}

\begin{corollary}\label{C:type152WURP}
Let $A$ be an algebra with almost permutable congruences. Then $\Conc A$
satisfies $\wurpe$.
\end{corollary}

\begin{remark}\label{Rk:type152WURP}
In case the distance $\delta$ satisfies the V-condition of type~$1$, the
statement $\wurpe$ in Theorem~\ref{T:type152WURP} can be strengthened
by taking $m=1$ in Definition~\ref{D:VWURPe}. Similarly, if $A$ is an
algebra with permutable congruences, then $\Conc A$ satisfies that
strengthening of $\wurpe$.
In particular, as any group, resp. any module, has permutable
congruences, both $\NSubc G$, for a group $G$, and $\Subc M$, for a
module~$M$, satisfy the strengthening of $\wurpe$ obtained by taking
$m=1$ in Definition~\ref{D:VWURPe}.
\end{remark}

As we shall see in Theorem~\ref{T:NonUnif}, not every distributive
\jzs\ can be join-generated by the range of a V-distance of type~$3/2$.
The situation changes dramatically for type~$2$. It is proved in \cite{GrLa} that any modular algebraic lattice is isomorphic to the congruence lattice
of an algebra with~$3$-permutable congruences. This easily implies the following result; nevertheless, we provide a much more direct argument, which will be useful for the proof of Theorem~\ref{T:ReprType2}.

\begin{proposition}\label{P:ReprType2}
Any distributive \jzs\ is the range of some V-distance of type~$2$.
\end{proposition}

\begin{proof}
Let $S$ be a distributive \jzs. We first observe that the map
$\mu_S\colon S\times S\to S$ defined by the rule
 \begin{equation}\label{Eq:DefmuS}
 \mu_S(x,y)=\begin{cases}
 x\vee y,&\text{ if }x\neq y,\\
 0,&\text{ if }x=y
 \end{cases}
 \end{equation}
is a surjective $S$-valued distance on $S$. Now suppose that we are
given a surjective $S$-valued distance $\delta\colon X\times X\to S$,
and let $x$, $y\in X$ and $\ba$, $\bb\in S$ such that
$\delta(x,y)\leq\ba\vee\bb$. Since $S$ is distributive, there are
$\ba'\leq\ba$ and $\bb'\leq\bb$ such that $\delta(x,y)=\ba'\vee\bb'$. We
put $X'=X\cup\set{u,v}$, where $u$ and $v$ are two distinct outside
points, and we extend $\delta$ to a distance $\delta'$ on $X'$ by
putting $\delta'(z,u)=\delta(z,x)\vee\ba'$ and
$\delta'(z,v)=\delta(z,y)\vee\ba'$, for all $z\in X$, while
$\delta'(u,v)=\bb'$. It is straightforward to verify that $\delta'$ is
an~$S$-valued distance on $X'$ extending $\delta$. Furthermore,
$\delta'(x,u)=\ba'\leq\ba$, $\delta'(u,v)=\bb'\leq\bb$, and
$\delta'(v,y)=\ba'\leq\ba$. Iterating this construction transfinitely,
taking direct limits at limit stages,
yields an~$S$-valued V-distance of type~$2$ extending~$\delta$.
\end{proof}

\section{Failure of $\wurpe$ in $\Conc F$, for $F$
free bounded lattice}\label{S:NonWURPE}

The main proof of the present section, that is, the proof of
Theorem~\ref{T:NonUnif}, follows the lines of the proofs of
\cite[Theorem~3.3]{PTW} and \cite[Corollary~2.1]{TuWe}. However, there are
a few necessary changes, mainly due to the new ``uniform refinement
property'' not being the same as the previously considered ones. As the
new result extends to any algebra, and not only lattices (see
Corollary~\ref{C:FreeNon}), we feel that it is still worthwhile to show
the main lines of the proof in some detail.

{}From now on until Lemma~\ref{L:eiUnique}, we shall fix
a non-distributive lattice variety $\cV$. For every set $X$, denote by
$\B_{\cV}(X)$ (or $\B(X)$ in case $\cV$ is understood) the
\emph{bounded} lattice in $\cV$ freely generated by chains $s_i<t_i$,
for $i\in X$. Note that if $Y$ is a subset of $X$, then there is
a unique retraction from $\B(X)$ onto $\B(Y)$, sending each~$s_i$ to~$0$
and each $t_i$ to $1$, for every $i\in X\setminus Y$. Thus, we shall
often identify $\B(Y)$ with the bounded sublattice of $\B(X)$ generated
by all $s_i$ and $t_i$ ($i\in Y$). Moreover, the abovementioned
retraction from
$\B(X)$ onto $\B(Y)$ induces a retraction from $\Conc\B(X)$
onto $\Conc\B(Y)$. Hence, we shall also identify $\Conc\B(Y)$
with the corresponding subsemilattice of $\Conc\B(X)$.

Now we fix a set $X$ such that $|X|\geq\aleph_2$. We denote, for all
$i\in X$, by $\ba_i$ and $\bb_i$ the compact congruences of $\B(X)$
defined by
\begin{equation}\label{Eq:Defnbaibi}
\ba_i=\Theta(0,s_i)\vee\Theta(t_i,1);\qquad
\bb_i=\Theta(s_i,t_i).
\end{equation}

In particular, note that $\ba_i\vee\bb_i=\one$, the largest congruence
of $\B(X)$.

Now, towards a contradiction, suppose that there are a positive integer
$n$, a decomposition $X=\bigcup\famm{X_k}{k<n}$, and a
family $\seqm{\bc_{i,j}}{\seq{i,j}\in X\times X}$ of elements of
$\Conc\B(X)$ witnessing the statement that $\Conc\B(X)$ satisfies
$\wurpe(\one)$, where~$\one$ denotes the largest congruence of $\B(X)$. We
pick $k<n$ such that $|X_k|=|X|$. By ``projecting everything on
$\B(X_k)$'' (as in \cite[page~224]{TuWe}), we might assume that
$X_k=X$.

Since the $\Conc$ functor preserves direct limits,
for all $i$, $j\in X$, there exists a finite subset $F(\set{i,j})$ of
$X$ such that both $\bc_{i,j}$ and $\bc_{j,i}$ belong to $\Conc\B(F(\set{i,j}))$.
By Kuratowski's Theorem, there are distinct elements $0$, $1$, $2$ of
$X$ such that $0\notin F(\set{1,2})$, $1\notin F(\set{0,2})$,
and $2\notin F(\set{0,1})$. Denote by
$\pi\colon\B(X)\twoheadrightarrow\B(\set{0,1,2})$ the canonical
retraction. For every $i\in\set{0,1,2}$, denote by $i'$ and
$i''$ the other two elements of~$\set{0,1,2}$, arranged in such a way
that $i'<i''$. We put $\bd_i=(\Conc\pi)(\bc_{i',i''})$, for all $i\in\set{0,1,2}$.

Applying the semilattice homomorphism $\Conc\pi$
to the inequalities satisfied by the elements $\bc_{i,j}$ yields
\begin{gather}
\bd_0\subseteq \ba_1\vee\ba_2,\bb_1\vee\bb_2;\qquad
\bd_1\subseteq\ba_0\vee\ba_2,\bb_0\vee\bb_2;\qquad
\bd_2\subseteq\ba_0\vee\ba_1,\bb_0\vee\bb_1;\label{Eq:D1}\\
\bd_0\vee\ba_2\vee\bb_1=
\bd_1\vee\ba_2\vee\bb_0=
\bd_2\vee\ba_1\vee\bb_0=\one;\label{Eq:D2}\\
\bd_1\subseteq\bd_0\vee\bd_2.\label{Eq:D3}
\end{gather}

As in \cite[Lemma~2.1]{PTW}, it is not hard to prove the following.

\begin{lemma}\label{L:Supp}
The congruence $\bd_i$ belongs to $\Conc\B(\set{i',i''})$, for all $i\in\set{0,1,2}$.
\end{lemma}

\begin{figure}[hbt]
\includegraphics{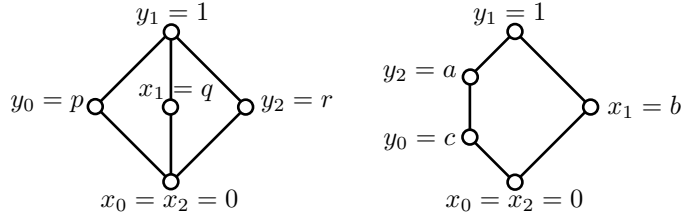}
\caption{The decorations of $M_3$ and $N_5$.}\label{Fig:DecM}
\end{figure}

Since $\cV$ is a non-distributive variety of lattices, it follows from
a classical result of lattice theory that $\cV$ contains as a member
some lattice $M\in\set{M_3,N_5}$. Decorate the lattice $M$
with three $2$-element chains $x_i<y_i$ (for $i\in\set{0,1,2}$) as in
\cite{PTW}, which we illustrate on Figure~\ref{Fig:DecM}.

The relevant properties of these decorations are summarized in
the two following straightforward lemmas.

\begin{lemma}\label{L:Ineq}
The decorations defined above satisfy the following
inequalities
 \begin{gather*}
 x_0\wedge y_1\leq x_1;\qquad y_1\leq x_1\vee y_0;\qquad
 x_1\wedge y_0\leq x_0;\qquad y_0\leq x_0\vee y_1;\\
 x_1\wedge y_2\leq x_2;\qquad y_2\leq x_2\vee y_1;\qquad
 x_2\wedge y_1\leq x_1;\qquad y_1\leq x_1\vee y_2,
 \end{gather*}
but $y_2\not\leq x_2\vee y_0$.
\end{lemma}

\begin{lemma}\label{L:ChDistr}
The sublattice of $M$ generated by $\set{x_{i'},x_{i''},y_{i'},y_{i''}}$ is distributive, for all $i\in\set{0,1,2}$.
\end{lemma}

\begin{figure}[htb]
\includegraphics{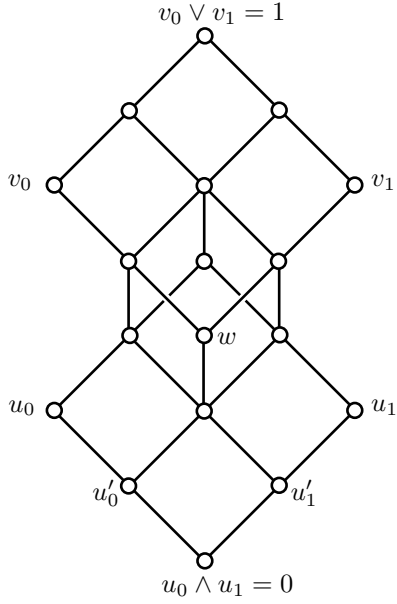}
\caption{The distributive lattice $D$.}\label{Fig:LattD}
\end{figure}

Now we shall denote by $D$ be the free product (i.e., the coproduct) of
two $2$-element chains, say $u_0<v_0$ and $u_1<v_1$, in the variety of
all \emph{distributive} lattices. The lattice $D$ is diagrammed on
Figure~\ref{Fig:LattD}.

The \jirr\
elements of $D$ are $u_0$, $u_1$, $v_0$, $v_1$, $u'_0=u_0\wedge v_1$,
$u'_1=u_1\wedge v_0$, and $w=v_0\wedge v_1$. Since $D$ is finite
distributive, its congruence lattice is finite Boolean, with seven atoms
$\boldsymbol{p}=\Theta_D(p_*,p)$, for $p\in\J(D)$ (where $p_*$ denotes
the unique lower cover of $p$ in $D$), that is,
 \begin{align*}
 \bu_0&=\Theta^+_D(u_0,v_1);&
 \bu_1&=\Theta^+_D(u_1,v_0);\\
 \bv_0&=\Theta^+_D(v_0,u_0\vee v_1);&
 \bv_1&=\Theta^+_D(v_1,u_1\vee v_0);\\
 \bu'_0&=\Theta^+_D(u_0\wedge v_1,u_1);&
 \bu'_1&=\Theta^+_D(u_1\wedge v_0,u_0);\\
 \bw&=\Theta_D((u_0\wedge v_1)\vee(u_1\wedge v_0),v_0\wedge v_1).
 \end{align*}
For all $i\in\set{0,1,2}$, let $\pi_i\colon\B(\set{i',i''})\to D$ be
the unique lattice homomorphism sending~$s_{i'}$ to~$u_0$,
$t_{i'}$ to~$v_0$, $s_{i''}$ to~$u_1$, $t_{i''}$ to~$v_1$. Furthermore,
denote by $\rho\colon\B(\set{0,1,2})\to M$ the unique lattice
homomorphism sending $s_i$ to~$x_i$ and $t_i$ to~$y_i$ (for all
$i\in\set{0,1,2}$); denote by $\rho_i$ the restriction of $\rho$ to
$\B(\set{i',i''})$.

We shall restate \cite[Lemma~3.1]{PTW} here for convenience.

\begin{lemma}\label{L:BasicDistr}
Let $L$ be any distributive lattice, let $a$, $b$, $a'$, $b'$
be elements of $L$. Then the equality
$\Theta^+_L(a,b)\cap\Theta^+_L(a',b')=\Theta^+_L(a\wedge a',b\vee b')$
holds.
\end{lemma}

Now we put $\be_i=(\Conc\pi_i)(\bd_i)$, for all $i\in\set{0,1,2}$.

\begin{lemma}\label{L:eiUnique}
The containments $\be^-\subseteq\be_i\subseteq\be^+$
hold for all $i\in\set{0,1,2}$, where we put
 \begin{align*}
 \be^-&=\Theta_D^+(u_0\wedge v_1,u_1)\vee\Theta_D^+(v_1,u_1\vee v_0),\\
 \be^+&=\Theta_D^+(u_0\wedge v_1,u_1)\vee\Theta_D^+(v_1,u_1\vee v_0)\vee
 \Theta_D^+(u_1\wedge v_0,u_0)\vee\Theta_D^+(v_0,u_0\vee v_1).
 \end{align*}
\end{lemma}

\begin{proof}
Applying $\Conc\pi_i$ to
the inequalities (\ref{Eq:D1}) and (\ref{Eq:D2}) yields
the following inequalities:
\begin{gather}
\be_i\subseteq\Theta(0,u_0)\vee\Theta(0,u_1)
\vee\Theta(v_0,1)\vee\Theta(v_1,1),\label{Eq:euv1}\\
\be_i\subseteq\Theta(u_0,v_0)\vee\Theta(u_1,v_1),\label{Eq:euv2}\\
\be_i\vee\Theta(0,u_1)\vee\Theta(v_1,1)\vee\Theta(u_0,v_0)=
\one.\label{Eq:euv3}
\end{gather}
By using Lemma~\ref{L:BasicDistr} and the distributivity of $\Con D$,
we obtain, by meeting \eqref{Eq:euv1} and \eqref{Eq:euv2}, the inequality
$\be_i\subseteq\be^+$. On the other hand, by using \eqref{Eq:euv3}
together with the equality
 \[
 \Theta(0,u_1)\vee\Theta(v_1,1)\vee\Theta(u_0,v_0)=
 \bu_0\vee\bu_1\vee\bu'_1\vee\bv_0\vee\bw,
 \]
(see Figure~\ref{Fig:LattD}), we obtain that
$\be^-=\bu'_0\vee\bv_1\subseteq\be_i$.
\end{proof}

Now, for all $i\in\set{0,1,2}$, it follows from Lemma~\ref{L:ChDistr} that
there exists a unique lattice homomorphism
$\varphi_i\colon D\to M$ such that
$\varphi_i\circ\pi_i=\rho_i$. Since $\Conc$ is a functor, we get from
this and from Lemma~\ref{L:eiUnique} that for all $i\in\set{0,1,2}$,
\begin{multline}
(\Conc\rho)(\bd_i)=(\Conc\varphi_i)(\be_i)\subseteq
(\Conc\varphi_i)(\be^+)\\
=(\Conc\varphi_i)
(\Theta^+(u_0\wedge v_1,u_1)\vee\Theta^+(v_1,u_1\vee v_0)\vee
\Theta^+(u_1\wedge v_0,u_0)\vee\Theta^+(v_0,u_0\vee v_1))\\
=\Theta^+(x_{i'}\wedge y_{i''},x_{i''})\vee
\Theta^+(y_{i''},x_{i''}\vee y_{i'})\vee
\Theta^+(x_{i''}\wedge y_{i'},x_{i'})\vee
\Theta^+(y_{i'},x_{i'}\vee y_{i''}).
\end{multline}
while
\begin{multline}
(\Conc\rho)(\bd_i)=(\Conc\varphi_i)(\be_i)\supseteq
(\Conc\varphi_i)(\be^-)\\
=(\Conc\varphi_i)
(\Theta^+(u_0\wedge v_1,u_1)\vee\Theta^+(v_1,u_1\vee v_0))\\
=\Theta^+(x_{i'}\wedge y_{i''},x_{i''})\vee
\Theta^+(y_{i''},x_{i''}\vee y_{i'}).
\end{multline}

In particular, we obtain, using Lemma~\ref{L:Ineq},
\begin{align*}
(\Conc\rho)(\bd_0)&=\zero,\\
(\Conc\rho)(\bd_2)&=\zero,\\
\text{while}\qquad(\Conc\rho)(\bd_1)&\supseteq
\Theta^+(x_0\wedge y_2,x_2)\vee\Theta^+(y_2,x_2\vee y_0)
\neq\zero.
\end{align*}
On the other hand, by applying $\Conc\rho$ to (\ref{Eq:D3}), we
obtain that
\[
(\Conc\rho)(\bd_1)\subseteq(\Conc\rho)(\bd_0)\vee(\Conc\rho)(\bd_2),
\]
a contradiction. Therefore, we have proved the following
theorem.

\begin{theorem}\label{T:NonUnif}
Let $\cV$ be any non-distributive variety of lattices, let $X$
be any set such that $|X|\geq\aleph_2$. Denote by $\B_{\cV}(X)$ the
free product in $\cV$ of $X$ copies of a two-element chain with
a least and a largest element added. Then $\Conc\B_{\cV}(X)$
does not satisfy $\wurpe$ at its largest element.
\end{theorem}

A ``local'' version of Theorem~\ref{T:NonUnif} is presented in
Theorem~\ref{T:NoLifttyp1}.

Observe that $\Conc\B_{\cV}(X)$, being the semilattice of compact
congruences of a lattice, is distributive.

As in \cite[Corollary~4.1]{PTW}, we obtain the following.

\begin{corollary}\label{C:FreeNon}
Let $L$ be any lattice that admits a lattice homomorphism
onto a free bounded lattice in the variety generated by either
$M_3$ or $N_5$ with $\aleph_2$ generators.
Then $\Conc L$ does not satisfy $\wurpe$. In particular, there exists no
V-distance of type~$3/2$ with range join-generating $\Conc L$. Hence
there is no algebra $A$ with almost permutable congruences such that
$\Con L\cong\Con A$.
\end{corollary}

\begin{proof}
The first part of the proof goes like the proof of
\cite[Corollary~4.1]{PTW}, using Lemma~\ref{L:WDUnif}. The rest of the
conclusion follows from Theorem~\ref{T:type152WURP}.
\end{proof}

\begin{corollary}\label{C:FreeNon2}
Let $\cV$ be any non-distributive variety of lattices
and let $F$ be any free \pup{resp., free bounded} lattice with at
least $\aleph_2$ generators in $\cV$. Then there exists no V-distance of
type~$3/2$ with range join-generating $\Conc F$. In particular, there is
no algebra~$A$ with almost permutable congruences such that
$\Con F\cong\Con A$.
\end{corollary}

By using Corollary~\ref{C:GrpMod}, we thus obtain the following.

\begin{corollary}\label{C:nonCGPiCMP}
Let $\cV$ be a non-distributive variety of lattices, let $F$ be any
free \pup{resp., free bounded} lattice with at least $\aleph_2$
generators in $\cV$, and put $D=\Con F$---a distributive, algebraic
lattice with $\aleph_2$ compact elements. Then there is no module~$M$
\pup{resp., no group~$G$, no $\ell$-group $G$} such that $\Sub M\cong D$
\pup{resp., $\NSub G\cong D$, $\Idl G\cong D$}.
\end{corollary}

Hence, not every distributive algebraic lattice is isomorphic to the
submodule lattice of some module, or to the normal subgroup lattice of
some group. However, our proof of this negative result requires at least
$\aleph_2$ compact elements. As we shall see in
Sections~\ref{S:PosCGMiCMP1} and~\ref{S:PosCGMiCMP2}, the $\aleph_2$ bound
is, in both cases of modules and groups, optimal.

\section{Representing distributive algebraic lattices with at most
$\aleph_1$ compact elements as submodule lattices of modules}
\label{S:PosCGMiCMP1}

In this section we deal with congruence lattices of right modules over
rings.

\begin{theorem}\label{T:SubMAl1}
Every distributive \jzs\ of size at most $\aleph_1$ is isomorphic to the
submodule lattice of some right module.
\end{theorem}

\begin{proof}
Let $S$ be a distributive \jzs\ of size at most $\aleph_1$. If $S$ has a
largest element, then it follows from the main result of \cite{Al1Reg}
that $S$ is isomorphic to the semilattice $\Idc R$ of all
finitely generated two-sided ideals of some (unital) von~Neumann regular
ring $R$.

In order to reduce ideals to submodules, we use a well-known trick.
As $R$ is a bimodule over itself, the tensor product
$\oll{R}=R^{\op}\otimes R$ can be endowed with a structure of (unital)
ring, with multiplication satisfying
$(a\otimes b)\cdot(a'\otimes b')=(a'a)\otimes(bb')$ (both $a'a$ and
$bb'$ are evaluated in $R$). Then $R$ is a right
$\oll{R}$-module, with scalar multiplication given by
$x\cdot(a\otimes b)=axb$, and the submodules of $R_{\oll{R}}$ are exactly
the two-sided ideals of $R$. Hence, $\Subc R_{\oll{R}}=\Idc R\cong S$.

In case $S$ has no unit, it is an ideal of the distributive
\jzus\ $S'=S\cup\set{1}$ for a new largest element $1$. By the previous
paragraph, $S'\cong\Subc M$ for some right module $M$, hence
$S\cong\Subc N$ where $N$ is the submodule of $M$ consisting of those
elements $x\in M$ such that the submodule generated by $x$ is sent
to an element of~$S$ by the isomorphism $\Subc M\cong S'$.
\end{proof}

The commutative case is quite different. For example, for a commutative von~Neumann regular ring~$R$, if $\Id R$ is finite, then, as it is distributive and complemented, it must be Boolean. In particular, the three-element chain is not isomorphic to the ideal lattice of any commutative von~Neumann regular ring. Even if regularity is removed, not every finite distributive lattice is allowed. For example, one can prove the following result: \emph{A finite distributive lattice~$D$ is isomorphic to the submodule lattice of a module over some commutative ring if{f}~$D$ is isomorphic to the ideal lattice of some commutative ring, if{f}~$D$ is
a product of chains}. In particular, the square $\two\times\two$ with a new bottom (resp., top) element added is not isomorphic to the submodule lattice of any module over a commutative ring.

\section{Representing distributive algebraic lattices with at most
$\aleph_1$ compact elements as normal subgroup lattices of groups}
\label{S:PosCGMiCMP2}

Every nonabelian simple group is ``neutral'' in the sense of \cite{FrMK}. Hence, the direction (1)$\Rightarrow$(5) in \cite[Theorem~8.5]{FrMK} yields the following well-known result, which holds despite the failure of congruence-distributivity in the variety of all groups.

\begin{lemma}\label{L:FinProdSimpl}
Let $n<\go$ and let $\seqm{G_i}{i<n}$ be a finite sequence of simple
nonabelian groups. Then the normal subgroups of $\prod_{i<n}G_i$ are
exactly the trivial ones, namely the products of the form
$\prod_{i<n}H_i$, where $H_i$ is either $G_i$ or $\set{1_{G_i}}$, for
all $i<n$. Consequently, $\NSub\left(\prod_{i<n}G_i\right)\cong\two^n$.
\end{lemma}

We denote by $\cF$ the class of all finite products of alternating
groups of the form $\mathfrak{A}_n$, for $n\geq5$. For a group
homomorphism $f\colon G\to H$, we denote by
$\NSub f\colon\NSub G\to\NSub H$ the \jzh\ that with any normal subgroup
$X$ of $G$ associates the normal subgroup of $H$ generated by
$f[X]$. The following square amalgamation result is crucial.
It is an analogue for groups of \cite[Theorem~1]{GrLW} (for lattices) or
\cite[Theorem~4.2]{Al1Reg} (for regular algebras over a division ring).

\begin{lemma}\label{L:SqAmalGrp}
Let $G_0$, $G_1$, $G_2$ be groups in $\cF$ and let
$f_1\colon G_0\to G_1$ and $f_2\colon G_0\to G_2$ be group
homomorphisms. Let $B$ be a finite Boolean semilattice, and, for
$i\in\set{1,2}$, let $\xg_i\colon\NSub G_i\to B$ be \jzh s such that
 \begin{equation}\label{Eq:xgifi}
 \xg_1\circ\NSub f_1=\xg_2\circ\NSub f_2.
 \end{equation}
Then there are a group $G$ in $\cF$, group homomorphisms
$g_i\colon G_i\to G$, for $i\in\set{1,2}$, and an isomorphism
$\alpha\colon\NSub G\to B$ such that $g_1\circ f_1=g_2\circ f_2$ and
$\alpha\circ\NSub g_i=\xg_i$ for all $i\in\set{1,2}$.
\end{lemma}

\begin{proof}[Outline of proof]
We follow the lines of the proofs of \cite[Theorem~1]{GrLW} or
\cite[Theorem~4.2]{Al1Reg}. First, by decomposing $B$ as a finite power
of $\two$, observing that~$\cF$ is closed under
finite direct products, and using Lemma~\ref{L:FinProdSimpl}, we reduce
to the case where $B=\two$, the two-element chain. Next, denoting by
$\xh$ the \jzh\ appearing on both sides of \eqref{Eq:xgifi}, we put
$G'_0=\setm{x\in G_0}{\xh([x])=0}$ (where $[x]$ denotes, again, the
normal subgroup generated by $x$), and, similarly, $G'_i=\setm{x\in
G_i}{\xg_i([x])=0}$, for $i\in\set{1,2}$. So
$G'_i$ is a normal subgroup of $G_i$, for all $i\in\set{0,1,2}$, and replacing $G_i$
by $G_i/G'_i$ makes it possible to reduce to the case where both $\xg_1$
and $\xg_2$ separate zero while both $f_1$ and $f_2$ are group
embeddings.

Hence the problem that we must solve is the following: given
group embeddings $f_i\colon G_0\into G_i$, for $i\in\set{1,2}$, we must
find a finite, simple, nonabelian group $G$ with group embeddings
$g_i\colon G_i\into G$ such that $g_1\circ f_1=g_2\circ f_2$. By the
positive solution of the amalgamation problem for finite groups (see
\cite[Section~15]{Neum54}), followed by embedding the resulting group
into some alternating group with index at least~$5$, this is possible.
\end{proof}

Now every distributive \jzs\ of size at most $\aleph_1$ is the direct
limit of some direct system of finite Boolean \jzs s and \jzh s;
furthermore, we may assume that the indexing set of the direct system is
a \emph{$2$-ladder}, that is, a lattice with zero where every interval
is finite and every element has at most two immediate predecessors.
Hence, by imitating the method of proof used in \cite[Theorem~2]{GrLW} or
\cite[Theorem~5.2]{Al1Reg}, it is not difficult to obtain the following
result.

\begin{theorem}\label{T:Al1GrpRepr}
Every distributive \jzs\ of size at most $\aleph_1$ is
isomorphic to the finitely generated normal subgroup semilattice of some
group which is a direct limit of members of $\cF$.
\end{theorem}

Reformulating the result in terms of algebraic lattices rather than
semilattices, together with the observation that all direct limits of
groups in $\cF$ are locally finite, gives the following.

\begin{corollary}\label{C:Al1GrpRepr}
Every distributive algebraic lattice with at most $\aleph_1$ compact
elements is isomorphic to the normal subgroup lattice of some locally
finite group.
\end{corollary}

\section{Representing distributive algebraic lattices with at most
$\aleph_0$ compact elements as $\ell$-ideal lattices of $\ell$-groups}
\label{S:PosCGMiCMP3}

The variety of $\ell$-groups is quite special, as it is both
congruence-distributive and congruence-permutable. The following lemma does not extend to the commutative case (for example, $\ZZ\times\ZZ$ cannot be embedded into any simple commutative $\ell$-group).

\begin{lemma}\label{L:EmbSLG}
Every $\ell$-group can be embedded into some simple $\ell$-group.
\end{lemma}

\begin{proof}
It follows from \cite[Corollary~5.2]{Pier72} that every $\ell$-group $G$
embeds into an $\ell$-group $H$ in which any two positive elements are
conjugate. In particular, $H$ is simple.
\end{proof}

The following result is a ``one-dimensional'' analogue for $\ell$-groups
of Lemma~\ref{L:SqAmalGrp}.

\begin{lemma}\label{L:1dimAmalLG}
For any $\ell$-group $G$, any finite Boolean semilattice $B$, and
any \jzh\ $\xf\colon\Idlc G\to B$, there are an $\ell$-group $H$, an
$\ell$-homomorphism $f\colon G\to\nobreak H$, and an isomorphism
$\alpha\colon\Idlc H\to B$ such that $\xf=\alpha\circ\Idlc f$.
\end{lemma}

\begin{proof}
Suppose first that $B=\two$. Observing that $I=\setm{x\in G}{\xf(G(x))=0}$
is an $\ell$-ideal of $G$, we let $H$ be any simple $\ell$-group
extending $G/I$ (see Lemma~\ref{L:EmbSLG}), we let $f\colon G\to H$ be the
composition of the canonical projection $G\onto G/I$ with the inclusion
map $G/I\into H$, and we let $\alpha\colon\Idlc H\to\two$ be the unique
isomorphism.

Now suppose that $B=\two^n$, for a natural number $n$. For each $i<n$, we
apply the result of the paragraph above to the $i$-th component
$\xf_i\colon\Idlc G\to\two$ of $\xf$, getting a simple $\ell$-group 
$H_i$, an $\ell$-homomorphism $f_i\colon G\to H_i$, and the isomorphism
$\alpha_i\colon\Idlc H_i\to\two$. Then we put $H=\prod_{i<n}H_i$,
$f\colon x\mapsto\seqm{f_i(x)}{i<n}$, and we let
$\alpha\colon\Idlc H\to\two^n$ be the canonical isomorphism.
\end{proof}

\begin{theorem}\label{T:Reprlgroups}
Every distributive at most countable \jzs\ is isomorphic to the
semilattice of all finitely generated $\ell$-ideals of some $\ell$-group.
\end{theorem}

Equivalently, every distributive algebraic lattice with (at most)
countably many compact elements is isomorphic to the $\ell$-ideal lattice
of some $\ell$-group.

\begin{proof}
It follows from \cite[Theorem~3.1]{BFMD} (see also
\cite[Theorem~6.6]{GW}) that every distributive at most countable \jzs\
$S$ can be expressed as the direct limit of a sequence
$\seqm{B_n}{n<\go}$ of finite Boolean semilattices, with all transition
maps $\xf_n\colon B_n\to B_{n+1}$ and limiting maps
$\xg_n\colon B_n\to S$ being \jzh s. We fix an $\ell$-group $G_0$ with an
isomorphism $\alpha_0\colon\Idlc G_0\onto B_0$. Suppose having
constructed an $\ell$-group $G_n$ with an isomorphism
$\alpha_n\colon\Idlc G_n\to B_n$. Applying Lemma~\ref{L:1dimAmalLG} to
$\xf_n\circ\alpha_n$, we obtain an $\ell$-group $G_{n+1}$,
an $\ell$-homomorphism $f_n\colon G_n\to G_{n+1}$, and an
isomorphism $\alpha_{n+1}\colon\Idlc G_{n+1}\to B_{n+1}$ such that
$\xf_n\circ\alpha_n=\alpha_{n+1}\circ\Idlc f_n$. Defining $G$ as the
direct limit of the sequence
 \[
\xymatrix{
G_0\ar[r]^{f_0} & G_1\ar[r]^{f_1} & G_2\ar[r]^(.4){f_2} & \dots\dots\ ,
}
 \]
an elementary categorical argument yields an isomorphism from $\Idlc G$
onto the direct limit $S$ of the sequence $\seqm{B_n}{n<\go}$.
\end{proof}

\section{Functorial representation by V-distances of type~$2$}
\label{S:Type2}

Observe that the argument of Proposition~\ref{P:ReprType2} is only a
small modification (with a more simple-minded proof) of B.~J\'onsson's
proof that every modular lattice has a type~$2$ representation, see
\cite{Jons53} or \cite[Theorem~IV.4.8]{GLT2}. It follows from
Corollary~\ref{C:FreeNon} that ``type~$2$'' cannot be improved to
``type~$1$''. In view of Proposition~\ref{P:Dist2Part}, this is somehow
surprising, as every distributive lattice has an embedding with
permutable congruences into some partition lattice. This illustrates the
observation that one can get much more from a distance than from an
embedding into a partition lattice.

We shall now present a strengthening of Proposition~\ref{P:ReprType2}
that shows that the construction can be made \emph{functorial}. We
introduce notations for the following categories:

\begin{itemize}
\item[(1)] $\cS$, the category of all distributive \jzs s with \jze s.

\item[(2)] $\cD$, the category of all surjective distances of the form
$\delta\colon X\times X\onto S$ with kernel the identity and $S$ a
distributive \jzs, with morphisms (see Definition~\ref{D:Vdist}) of the
form $\seq{f,\xf}\colon\seq{X,\lambda}\to\seq{Y,\mu}$ with both $f$
and $\xf$ one-to-one.

\item[(3)] $\cD_2$, the full subcategory of $\cD$ consisting of all
V-distances of type~$2$.
\end{itemize}

Furthermore, denote by $\Pi\colon\cD\to\cS$ the forgetful functor
(see Definition~\ref{D:Vdist}).

\begin{theorem}\label{T:ReprType2}
There exists a direct limits preserving functor $\Phi\colon\cS\to\cD_2$
such that the composition $\Pi\circ\Phi$ is equivalent to the identity.
\end{theorem}

Hence the functor $\Phi$ assigns to each distributive \jzs\ $S$ a set
$X_S$ together with a surjective $S$-valued V-distance
$\delta_S\colon X_S\times X_S\onto S$ of type~$2$.

\begin{proof}
The proof of Proposition~\ref{P:ReprType2} depends of the enumeration
order of a certain transfinite sequence of quadruples $\seq{x,y,\ba,\bb}$,
which prevents it from being functorial. We fix this by adjoining all
such quadruples simultaneously, and by describing the corresponding
extension. So, for a distance $\delta\colon X\times X\to S$, we put
$S^-=S\setminus\set{0}$, and
 \[
 \cH(\delta)=\setm{\seq{x,y,\ba,\bb}\in X\times X\times S^-\times S^-}
 {\delta(x,y)=\ba\vee\bb}.
 \]
For $\xi=\seq{x,y,\ba,\bb}\in\cH(\delta)$, we put $x_{\xi}^0=x$,
$x_{\xi}^1=y$, $\ba_{\xi}=\ba$, and $\bb_{\xi}=\bb$.
Now we put
$X'=X\cup\setm{u_{\xi}^i}{\xi\in\cH(\delta)\text{ and }i\in\set{0,1}}$,
where the elements $u_{\xi}^i$ are pairwise distinct symbols outside
$X$. We define a map $\delta'\colon X'\times X'\to S$ by requiring
$\delta'$ to extend $\delta$, with value zero on the diagonal, and by the
rule
 \begin{align*}
 \delta'(u_{\xi}^i,u_{\eta}^j)&=\begin{cases}
 |i-j|\cdot\bb_{\xi},&\text{ if }\xi=\eta,\\
 \ba_{\xi}\vee\ba_{\eta}\vee\delta(x_{\xi}^i,x_{\eta}^j),
 &\text{ if }\xi\neq\eta,
 \end{cases}\\
 \delta'(u_{\xi}^i,z)
 =\delta'(z,u_{\xi}^i)&=\delta(z,x_{\xi}^i)\vee\ba_{\xi},
 \end{align*}
for all $\xi$, $\eta\in\cH(\delta)$, all $i$, $j\in\set{0,1}$, and all
$z\in X$.

It is straightforward, though somewhat tedious, to verify that
$\delta'$ is an~$S$-valued distance on $X'$, that it extends~$\delta$, and
that its kernel is the identity of $X'$ in case the kernel of $\delta$ is
the identity of $X$ (because the semilattice elements $\ba_{\xi}$ and
$\bb_{\xi}$ are nonzero). Furthermore, if $S$ is distributive, then every
V-condition problem for~$\delta$ of the form
$\delta(x,y)\leq\ba\vee\bb$ can be refined to a problem of the form
$\delta(x,y)=\ba'\vee\bb'$, for some $\ba'\leq\ba$ and $\bb'\leq\bb$
(because $S$ is distributive), and such a problem has a solution of
type~$2$ for $\delta'$. Namely, in case both
$\ba'$ and $\bb'$ are nonzero (otherwise the problem can be solved in
$X$), put $\xi=\seq{x,y,\ba',\bb'}$, and observe that
$\delta'(x,u_{\xi}^0)=\ba'$, $\delta'(u_{\xi}^0,u_{\xi}^1)=\bb'$, and
$\delta'(u_{\xi}^1,y)=\ba'$.

Hence, if we put $\seq{X_0,\delta_0}=\seq{X,\delta}$, then
$\seq{X_{n+1},\delta_{n+1}}=\seq{(X_n)',(\delta_n)'}$ for all $n<\go$,
and finally $\ol{X}=\bigcup\famm{X_n}{n<\go}$ and
$\ol{\delta}=\bigcup\famm{\delta_n}{n<\go}$, the pair
$\Psi(\seq{X,\delta})=\seq{\ol{X},\ol{\delta}}$ is an~$S$-valued
V-distance of type~$2$ extending $\seq{X,\delta}$. Every morphism
$\seq{f,\xf}\colon\seq{X,\lambda}\to\seq{Y,\mu}$ in $\cS$ extends
canonically to a morphism
$\seq{f',\xf}\colon\seq{X',\lambda'}\to\seq{Y',\mu'}$ (the underlying
semilattice map $\xf$ is the same), by defining
 \[
 f'(u_{\xi}^i)=u_{f\xi}^i,\text{ for all }\xi\in\cH(\lambda)
 \text{ and all }i<2,
 \]
where we put, of course,
 \[
 f\seq{x,y,\ba,\bb}=\seq{f(x),f(y),\xf(\ba),\xf(\bb)},\quad\text{for all }
 \seq{x,y,\ba,\bb}\in\cH(\lambda).
 \]
Hence, by an easy induction argument, $\seq{f,\xf}$ extends canonically
to a morphism $\Psi(\seq{f,\xf})=\seq{\oll{f},\xf}\colon
\seq{\oll{X},\oll{\lambda}}\to\seq{\oll{Y},\oll{\mu}}$, and the
correspondence $\seq{f,\xf}\mapsto\seq{\oll{f},\xf}$ is itself a
functor. As the construction defining the correspondence
$\seq{X,\delta}\mapsto\seq{X',\delta'}$ is local, the functor~$\Psi$
preserves direct limits.

It remains to find something to start with, to which we can apply
$\Psi$. A possibility is to use the distance $\mu_S$, given by
\eqref{Eq:DefmuS}, introduced in the proof
of Proposition~\ref{P:ReprType2}. The correspondence
$S\mapsto\mu_S$ defines a functor, in particular, if $f\colon S\into T$
is an embedding of distributive
\jzs s, then the equality $\mu_T(f(x),f(y))=f(\mu_S(x,y))$ holds, for
all $x$, $y\in S$. The desired functor $\Phi$ is given by
$\Phi(S)=\Psi(\seq{S,\mu_S})$, for any distributive \jzs\ $S$.
\end{proof}

In contrast with the result of Theorem~\ref{T:ReprType2}, we shall
isolate a finite, ``combinatorial'' reason for the forgetful functor from
V-distances of type~$3/2$ to distributive \jzs s not to admit any left
inverse. By contrast, we recall that for V-distances of type~$2$, the
corresponding result is positive, see Theorem~\ref{T:ReprType2}. In order
to establish the negative result, we shall use the example $\Dac$ of
\cite[Section~7]{TuWe}, and extend the corresponding result from lattices
with almost permutable congruences to arbitrary V-distances of type~$3/2$.

\begin{figure}[htb]
 \[
{
\def\labelstyle{\displaystyle}
\xymatrixrowsep{2pc}\xymatrixcolsep{1.5pc}
\xymatrix{
& \Pow(8) & \\
\Pow(4)\ar[ru]^{\xh_2} & \Pow(4)\ar[u]|-{\xh_1}&
\Pow(4)\ar[lu]_{\xh_0}\\
& &\\
\Pow(2)\ar[uu]^{\xf}\ar[ruu]|-(.3){\strut\xf}
& \Pow(2)\ar[luu]|-(.7){\strut\xg}\ar[ruu]|-(.7){\strut\xf} &
\Pow(2)\ar[uu]_{\xg}\ar[luu]|-(.3){\strut\xg}\\
&\two\ar[lu]^{\be}\ar[u]|-{\strut\be}\ar[ru]\ar[ru]_{\be}& }
}
 \]
\caption{The cube $\Dac$, unliftable by V-distances of type~$3/2$.}
\label{Fig:UnliftCube}
\end{figure}
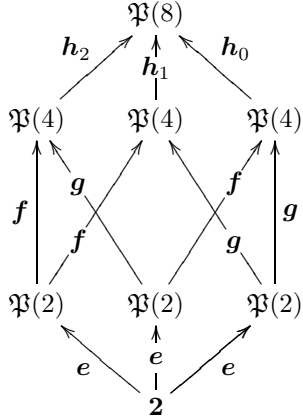

We recall that $\Dac$ is the (commutative) cube of finite Boolean
semilattices represented on Figure~\ref{Fig:UnliftCube}, where $\Pow(X)$
denotes the powerset algebra of a set $X$ and $\be$, $\xf$, $\xg$,
$\xh_0$, $\xh_1$, and $\xh_2$ are the \jzh s (and, in fact, \jzue s)
defined by their values on atoms as follows:
 \begin{align*}
 &\be(1)=\set{0,1},&&\\
 \xf\colon\begin{cases}
 \set{0}\mapsto\set{0,1},\\ \set{1}\mapsto\set{2,3},
 \end{cases}&
 \xg\colon\begin{cases}
 \set{0}\mapsto\set{0,2},\\ \set{1}\mapsto\set{1,3},
 \end{cases}\\
 \xh_0\colon\begin{cases}
 \set{0}\mapsto\set{0,4,7},\\
 \set{1}\mapsto\set{3,5,6},\\
 \set{2}\mapsto\set{2,5,6},\\
 \set{3}\mapsto\set{1,4,7},
 \end{cases}&
  \xh_1\colon\begin{cases}
 \set{0}\mapsto\set{0,4,5,7},\\
 \set{1}\mapsto\set{1,4,6,7},\\
 \set{2}\mapsto\set{2,5,6,7},\\
 \set{3}\mapsto\set{3,4,5,6},
 \end{cases}&
 \xh_2\colon\begin{cases}
 \set{0}\mapsto\set{0,4,6},\\
 \set{1}\mapsto\set{1,5,7},\\
 \set{2}\mapsto\set{3,5,7},\\
 \set{3}\mapsto\set{2,4,6},
 \end{cases}&
 \end{align*}

\begin{theorem}\label{T:NoLifttyp1}
The diagram $\Dac$ has no lifting, with
respect to the forgetful functor, by distances, surjective at level $0$
and satisfying the V-condition of type~$3/2$ at level~$1$.
\end{theorem}

\begin{proof}
Suppose that the diagram of Figure~\ref{Fig:UnliftCube} is lifted by a
diagram of distances, with distances $\lambda\colon X\times X\to\two$,
$\lambda_i\colon X_i\times X_i\to\Pow(2)$,
$\mu_i\colon Y_i\times Y_i\to\Pow(4)$, and
$\mu\colon Y\times\nobreak Y\to\nobreak\Pow(8)$,
for all $i\in\set{0,1,2}$, see Figure~\ref{Fig:DistDiagr}.

\begin{figure}[htb]
 \[
{
\def\labelstyle{\displaystyle}
\xymatrixrowsep{2pc}\xymatrixcolsep{1.5pc}
\xymatrix{
& \mu & \\
\mu_2\ar[ru] & \mu_1\ar[u]& \mu_0\ar[lu]\\
\lambda_0\ar[u]\ar[ru] & \lambda_1\ar[lu]\ar[ru] &
\lambda_2\ar[u]\ar[lu]\\
&\lambda\ar[lu]\ar[u]\ar[ru]&
}
}
 \]
\caption{A commutative diagram of
semilattice-valued distances.}\label{Fig:DistDiagr}
\end{figure}
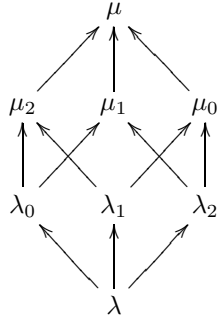

We assume that $\lambda$ is surjective and that
$\lambda_i$ is a V-distance of type~$3/2$, for all $i\in\set{0,1,2}$.
Denote by $f_{U,V}$ the canonical map from $U$ to $V$
given by this lifting, for $U$ below~$V$ among $X$, $X_0$, $X_1$, $X_2$,
$Y_0$, $Y_1$, $Y_2$, $Y$. After having replaced each of those sets $U$ by
its quotient by the kernel of the corresponding distance, and then by its
image in $Y$ under $f_{U,Y}$, we may assume that
$f_{U,V}$ is the inclusion map from~$U$ into~$V$, for all~$U$ below $V$
among $X$, $X_0$, $X_1$, $X_2$, $Y_0$, $Y_1$, $Y_2$, $Y$.

Since $\lambda$ is surjective, there are $x$, $y\in X$ such that
$\lambda(x,y)=1$. For all $i\in\set{0,1,2}$,
 \[
\lambda_i(x,y)=\be(\lambda(x,y))=\be(1)=\set{0,1}=\set{0}\cup\set{1},
 \]
thus, since $\lambda_i$ satisfies the V-condition of type~$3/2$, there
exists $z_i\in X_i$ such that
 \begin{equation}\label{Eq:PQ(i)}
 \begin{aligned}
 \text{either }&\lambda_i(x,z_i)=\set{0}\text{ and }
 \lambda_i(z_i,y)=\set{1}\quad(\text{say, }P(i))\\
 \text{or }&\lambda_i(x,z_i)=\set{1}\text{ and }
 \lambda_i(z_i,y)=\set{0}\quad(\text{say, }Q(i)).
 \end{aligned}
 \end{equation}
So we have eight cases to consider, according to which combination of $P$
and $Q$ occurs in \eqref{Eq:PQ(i)} for $i\in\set{0,1,2}$. In each case, we shall
obtain the inequality
 \begin{equation}\label{Eq:NegTr}
 \mu(z_0,z_2)\not\subseteq\mu(z_0,z_1)\cup\mu(z_1,z_2),
 \end{equation}
which will contradict the triangular inequality for $\mu$.

\Case{1} $P(0)$, $P(1)$, and $P(2)$ hold.
Then $\mu_2(z_0,x)=\xf\lambda_0(x,z_0)=\set{0,1}$ and
$\mu_2(x,z_1)=\xg(\lambda_1(x,z_1))=\set{0,2}$, whence
$\mu_2(z_0,z_1)\subseteq\set{0,1,2}$. Similarly, replacing~$x$ by~$y$ in
the argument above, $\mu_2(z_0,y)=\xf\lambda_0(z_0,y)=\set{2,3}$ and
$\mu_2(y,z_1)=\xg(\lambda_1(z_1,y))=\set{1,3}$, whence
$\mu_2(z_0,z_1)\subseteq\set{1,2,3}$. Therefore,
$\mu_2(z_0,z_1)\subseteq\set{1,2}$. On the other hand, from
$\mu_2(x,z_0)\cup\mu_2(z_0,z_1)=\mu_2(x,z_1)\cup\mu_2(z_0,z_1)$ the
converse inclusion follows, whence $\mu_2(z_0,z_1)=\set{1,2}$. Similar
computations yield that $\mu_1(z_0,z_2)=\mu_0(z_1,z_2)=\set{1,2}$.

Hence, we obtain the equalities
 \begin{align*}
 \mu(z_0,z_1)&=\xh_2\mu_2(z_0,z_1)=\set{1,3,5,7},\\
 \mu(z_0,z_2)&=\xh_1\mu_1(z_0,z_2)=\set{1,2,4,5,6,7},\\
 \mu(z_1,z_2)&=\xh_0\mu_0(z_1,z_2)=\set{2,3,5,6}.
 \end{align*}
Observe that $4$ belongs to $\mu(z_0,z_2)$ but not to
$\mu(z_0,z_1)\cup\mu(z_1,z_2)$.

\Case{2} $P(0)$, $P(1)$, and $Q(2)$ hold. As in Case~1, we obtain
 \[
 \mu_2(z_0,z_1)=\set{1,2}\text{ and }
 \mu_1(z_0,z_2)=\mu_0(z_1,z_2)=\set{0,3},
 \]
thus $\mu(z_0,z_1)=\set{1,3,5,7}$, $\mu(z_0,z_2)=\set{0,3,4,5,6,7}$,
and $\mu(z_1,z_2)=\set{0,1,4,7}$, which confirms \eqref{Eq:NegTr} and
thus causes a contradiction.

\Case{3} $P(0)$, $Q(1)$, and $P(2)$ hold. We obtain
 \[
 \mu_2(z_0,z_1)=\mu_0(z_1,z_2)=\set{0,3}\text{ and }
 \mu_1(z_0,z_2)=\set{1,2},
 \]
thus $\mu(z_0,z_1)=\set{0,2,4,6}$, $\mu(z_0,z_2)=\set{1,2,4,5,6,7}$,
and $\mu(z_1,z_2)=\set{0,1,4,7}$.

\Case{4} $P(0)$, $Q(1)$, and $Q(2)$ hold. We obtain
 \[
 \mu_2(z_0,z_1)=\mu_1(z_0,z_2)=\set{0,3}\text{ and }
 \mu_0(z_1,z_2)=\set{1,2},
 \]
thus $\mu(z_0,z_1)=\set{0,2,4,6}$, $\mu(z_0,z_2)=\set{0,3,4,5,6,7}$,
and $\mu(z_1,z_2)=\set{2,3,5,6}$.

\Case{5} $Q(0)$, $P(1)$, and $P(2)$ hold. We obtain
 \[
 \mu_2(z_0,z_1)=\mu_1(z_0,z_2)=\set{0,3}\text{ and }
 \mu_0(z_1,z_2)=\set{1,2},
 \]
thus $\mu(z_0,z_1)=\set{0,2,4,6}$, $\mu(z_0,z_2)=\set{0,3,4,5,6,7}$,
and $\mu(z_1,z_2)=\set{2,3,5,6}$.

\Case{6} $Q(0)$, $P(1)$, and $Q(2)$ hold. We obtain
 \[
 \mu_2(z_0,z_1)=\mu_0(z_1,z_2)=\set{0,3}\text{ and }
 \mu_1(z_0,z_2)=\set{1,2},
 \]
thus $\mu(z_0,z_1)=\set{0,2,4,6}$, $\mu(z_0,z_2)=\set{1,2,4,5,6,7}$,
and $\mu(z_1,z_2)=\set{0,1,4,7}$.

\Case{7} $Q(0)$, $Q(1)$, and $P(2)$ hold. We obtain
 \[
 \mu_2(z_0,z_1)=\set{1,2}\text{ and }
 \mu_1(z_0,z_2)=\mu_0(z_1,z_2)=\set{0,3},
 \]
thus $\mu(z_0,z_1)=\set{1,3,5,7}$, $\mu(z_0,z_2)=\set{0,3,4,5,6,7}$,
and $\mu(z_1,z_2)=\set{0,1,4,7}$.

\Case{8} $Q(0)$, $Q(1)$, and $Q(2)$ hold. We obtain
 \[
 \mu_2(z_0,z_1)=\mu_1(z_0,z_2)=\mu_0(z_1,z_2)=\set{1,2},
 \]
thus $\mu(z_0,z_1)=\set{1,3,5,7}$, $\mu(z_0,z_2)=\set{1,2,4,5,6,7}$,
and $\mu(z_1,z_2)=\set{2,3,5,6}$.

In all cases, we obtain a contradiction.
\end{proof}

A ``global'' version of Theorem~\ref{T:NoLifttyp1} is presented in
Theorem~\ref{T:NonUnif}.

The following corollary extends \cite[Theorem~7.1]{TuWe} from lattices to
arbitrary algebras.

\begin{corollary}\label{C:NoLifttyp1}
The diagram $\Dac$ has no lifting, with respect to the congruence lattice
functor, by algebras with almost permutable congruences.
\end{corollary}

About other commonly encountered structures, we obtain the following.

\begin{corollary}\label{C:NoLifttyp1'}
The diagram $\Dac$ has no lifting by groups with respect to the $\NSub$
functor, and no lifting by modules with respect to the $\Sub$ functor.
\end{corollary}

The following example offers a significant difference between the
situations for groups and modules.

\begin{example}\label{Ex:NoLift1}
The diagonal map $\two\into\two^2$ has no lifting, with respect to the
$\Sub$ functor, by modules over any ring. Indeed, suppose that $A\into
B\times C$ is such a lifting, with $A$, $B$, and $C$ simple modules.
Projecting on $B$ and on $C$ yields that~$A$ is isomorphic to a submodule
of both $B$ and $C$, whence, by simplicity, $A$, $B$, and $C$ are pairwise
isomorphic. But then, $B\times C\cong B\times B$ has the diagonal as a
submodule, so its submodule lattice cannot be isomorphic to $\two^2$.

By contrast, every square of finite Boolean \jzs s can be lifted, with
respect to the $\NSub$ functor, by \emph{groups}, see
Lemma~\ref{L:SqAmalGrp}.
\end{example}

\section{Open problems}\label{S:Pbs}

Although we do know that the negative result of
Corollary~\ref{C:FreeNon2} applies to $\ell$-groups (for every
$\ell$-group has permutable congruences), we do not know
whether the positive results proved here for modules
(Theorem~\ref{T:SubMAl1}) or for groups (Theorem~\ref{T:Al1GrpRepr})
extend to $\ell$-groups. The problem is that the class of all
$\ell$-groups does not satisfy the amalgamation property, see
\cite[Theorem~3.1]{Pier72}, so the proof of Lemma~\ref{L:SqAmalGrp}
cannot be used in this context, and so we do not know how to extend
Theorem~\ref{T:Reprlgroups} to the first uncountable level.

\begin{problem}\label{Pb:lGroups}
Is every distributive algebraic lattice with $\aleph_1$ compact
elements isomorphic to the $\ell$-ideal lattice of some
$\ell$-group?
\end{problem}

Our next question is related to the functor $\Phi$ obtained in the
statement of Theorem~\ref{T:ReprType2}.

\begin{problem}\label{Pb:FinTyp2}
Does there exist a functor $\Phi$ as in Theorem~\ref{T:ReprType2} that
sends finite semilattices to distances with finite underlying sets?
\end{problem}

That is, can we assign \emph{functorially} (with respect to \jze s), to
each finite distributive \jzs\ $S$, a surjective V-distance
$\seq{X_S,\delta_S}$ of type~$2$ with
$\delta_S\colon X_S\times X_S\onto S$ and $X_S$ \emph{finite}?

\end{document}